\documentclass[12pt]{amsart}

\usepackage{latexsym}
\usepackage{a4}
\usepackage{amssymb}
\usepackage{amsmath}
\usepackage[english]{babel}
\usepackage{tikz}
\usepackage{bm}
\usepackage{hyperref}

\newtheorem{thm}{Theorem}[section]
\newtheorem{pro}[thm]{Proposition}
\newtheorem{lem}[thm]{Lemma}

\newtheorem{cor}[thm]{Corollary}
\theoremstyle{definition}
\newtheorem{exa}[thm]{Example}
\newtheorem{de}[thm]{Definition}
\newtheorem{rem}[thm]{Remark}
\numberwithin{equation}{section}

\newcommand{\N}{\mathbb{N}}

\newcommand{\Con}{\operatorname{Con}}

\newcommand{\Pol}{\operatorname{Pol}}

\newcommand{\ar}{\operatorname{ar}}
\newcommand{\bS}{\mathbf{S}}
\newcommand{\bA}{\mathbf{A}}
\newcommand{\bG}{\mathbf{G}}

\newcommand{\oa}{\mathbf{a}}
\newcommand{\ob}{\mathbf{b}}
\newcommand{\oc}{\mathbf{c}}
\newcommand{\od}{\mathbf{d}}

\newcommand{\ox}{\mathbf{x}}

\newcommand{\oy}{\mathbf{y}}
\newcommand{\oz}{\mathbf{z}}

\title{Supernilpotent semigroups}

\thanks{Supported by the Ministry of Science, Education and Technological Development of the Republic of Serbia (Grant No. 451-03-68/2022-14/200125) and the Ministry of Scientific and Technological Development, Higher Education and Information Society of the Republic of Srpska (Grant No. 19.041/68-13-27/2023.)}
\date{\today}

\begin{document}
\bibliographystyle{amsalpha}

\maketitle

\begin{center}
\begin{small}
Jelena Radovi\'{c}\\
Department of Mathematics\\
University of East Sarajevo\\
71123 East Sarajevo\\
Bosnia and Herzegovina\\
ORCID: 0000-0002-1023-4463\\
{\tt jelena.radovic@ff.ues.rs.ba}\\[5mm]
\end{small}
\end{center}
\begin{center}
\begin{small}
Neboj\v{s}a Mudrinski\footnote{Corresponding author}\\
Department of Mathematics and Informatics\\
Faculty of Sciences\\
University of Novi Sad\\
21000 Novi Sad\\
Serbia\\
ORCID: 0000-0001-9830-6603\\
{\tt nmudrinski@dmi.uns.ac.rs}\\[5mm]
\end{small}
\end{center}

\begin{abstract}
Regular abelian semigroups are isomorphic to a direct product of an abelian group and a rectangular band (Warne, 1994). Seeking for a similar result for nilpotency, solvability and supernilpotency of regular semigroups, we obtain that analogous statement is true only in orthodox semigroups. We provide an example that shows that the same does not have to be true in regular semigroups that are not orthodox.
\end{abstract}

\maketitle

\section{Motivation}

In this note we study properties of semigroups that are defined by the commutator, see  (Definition \ref{DefnHC}). Namely, we are interested when is the commutator or its iteration of the full congruence 1 equal to the equality relation 0. The first such property is the \emph{term condition} or \emph{abelian}, that is,
when $[1,1]=0$. Abelian semigroups have been completely characterized by Warne in \cite{Warne}. In particular, as in Proposition \ref{WarneCor2.6}, a regular semigroup $\bS$ is abelian if and only if it is isomorphic to a direct product $\bG\times \mathbf{A}\times \mathbf{B}$, where $\bG$ is an abelian group, $\mathbf{A}$ is a left zero semigroup (multiplication is the first projection) and $\mathbf{B}$ is a right zero semigroup (multiplication is the second projection). Further properties well developed in group theory, and studied in algebras in general, that are defined using iterated commutator of the full congruence $1$ are \textit{nilpotency} and
\textit{solvability} (Definition \ref{DefNilpSolv}). These properties can be further generalized trough the notion of
\textit{supernilpotency} (Definition \ref{supernilpotent}). It is important to note that in general case, a supernilpotent algebra does not have to be nilpotent (\cite{MM:SNNIN}). Nilpotency, solvability and supernilpotency are weaker conditions then abelian, see Remark \ref{remNilp} and Proposition \ref{HC3}. One can ask whether
a characterization of Warne type can hold for such properties. Namely, we wonder whether the word abelian in Warne's
characterization for regular semigroups can be replaced by the word nilpotent or solvable or supernilpotent. If this is true, more precisely if the regular semigroup is isomorphic to the direct product of a group, left-zero semigroup and right zero-semigroup, then it should be an orthodox semigroup, see Proposition \ref{PropBandTimesGroup}. Orthodox semigroups are introduced as in Definition \ref{DefnOrt}. Unfortunately, neither nilpotency nor supernilpotency in regular
semigroups does not force it to be orthodox, see Example \ref{nilpotentnotorthodox}. Consequently, using Remark \ref{remNilp} we know that solvability does not force regular semigroup to be orthodox either. Hence if we seek for Warne's type result in regular semigroups for the mentioned weaker conditions, we should assume just orthodox semigroups. The main results of this note are that orthodox semigroup is $n$-nilpotent ($n$-solvable, $n$-supernilpotent) if and only if it is isomorphic to a direct product $\bG\times \mathbf{A}\times \mathbf{B}$ where
$\bG$ is an $n$-nilpotent ($n$-solvable, $n$-nilpotent) group, $\mathbf{A}$ is a left zero semigroup and $\mathbf{B}$ is a right zero semigroup (Proposition \ref{ThmOrthodoxNilpSolv} and Theorem \ref{SupernilpOrthodox}). Therefore, we obtain that in orthodox semigroups, conditions of nilpotency and supernilpotency are equivalent. As a corollary we have that the same is true in inverse semigroups, see Proposition \ref{inverse}.

\section{Preliminaries}

In this section we will review some definitions and properties of regular semigroups and commutator theory, that are needed for our results. For more details one can refer to monographs \cite{Howie} and \cite{McMc}.

A semigroup $\bS=(S,\cdot)$ is \emph{regular} if for every element $x$ there exists an element $y\in S$ such that $xyx=x$. For an element $x\in S$ we call a $y\in S$ such that $xyx=x$ and $yxy=y$ an \emph{inverse} of $x$. In this note we deal with the following subclasses of regular semigroups: rectangular bands, completely simple semigroups and orthodox semigroups. A semigroup $\bS$ is a \emph{rectangular band} if and only if it is isomorphic to a direct product of a left zero semigroup and a right zero semigroup, that is, to a semigroup of the form $L\times R$ where $L$, $R$ are nonempty sets, and where multiplication is given by $(l_1,r_1)\cdot (l_2,r_2)=(l_1,r_2)$ (\cite[Theorem 1.1.3]{Howie}).

Let us note that in semigroups there is no one-to-one correspondence between congruences and ideals, more precisely not every congruence is determined by an ideal. Semigroups without proper ideals are called \emph{simple}. An \emph{idempotent} in $S$ is every $a\in S$ such that $a^2=a$. The set of all idempotents of the semigroup $\bS$ will be, as usual, denoted by $E=E(S)$. Let us recall that there is a \textit{natural partial order} on the set of idempotents $E(S)$ of semigroup $\bS$, defined by $e\leq f$ if and only if $ef=fe=e$.  In a semigroup without zero, an idempotent $e$ is \textit{primitive} if it is minimal in the natural ordering of idempotents. If a simple semigroup $\bS$ contains a primitive idempotent, we say that $\bS$ is a \textit{completely simple semigroup}. Completely simple semigroups are well studied, and characterized by the following statement.

\begin{pro}\textup{(\cite{Howie} Theorem 3.3.1, Theorem 3.4.2)}\label{HowieTe3.3.1}
Let $G$ be a group with the identity element $e$, let $I$ and $\Lambda$ be non-empty sets such that there exists an element $1\in I\cap\Lambda$, and let $P=[p_{\lambda i}]$ be a $\Lambda\times I$ matrix with entries in $G$, that satisfies the property $p_{1,i}=e=p_{\lambda,1}$ for every $i\in I,\lambda\in\Lambda$. Let $S=I\times G\times\Lambda$, and define multiplication on $S$ by
\begin{equation*}
(i,g,\lambda) \cdot (j,h,\mu) = (i,gp_{\lambda j} h,\mu).
\end{equation*}
Then $\mathcal{M}[G;I,\Lambda;P]=(S,\cdot)$ is a completely simple semigroup. Conversely, every completely simple semigroup is isomorphic to a semigroup constructed in this way.
\end{pro}

Semigroup $\mathcal{M}[G;I,\Lambda;P]$ described in previous proposition is usually called $I\times\Lambda$ \emph{Rees matrix semigroup over the group $G$, with normal structure matrix $P$}, see \cite[Chapters 3.3, 3.4]{Howie}.

We recall that every congruence $\rho$ on a completely simple semigroup $\bS=\mathcal{M}[G;I,\Lambda;P]$ determines equivalence relations $\rho_I$ and $\rho_{\Lambda}$ on sets $I$ and $\Lambda$, respectively, in the following way
\begin{align*}
    & \rho_I = \{ (i,j)\in I\times I : (i,e,1)\, \rho \, (j,e,1)\}, \\
    & \rho_{\Lambda} = \{ (\lambda,\mu) \in \Lambda\times\Lambda : (1,e,\lambda)\, \rho \, (1,e,\mu) \}.
\end{align*}
A congruence $\rho$ of $\bS$ determines the linked triple $(\rho_I,N_{\rho},\rho_{\Lambda})$ where $N_{\rho}=\{g\in
G\,:\,(1,g,1)\,\rho\,(1,e,1)\}$ is the normal subgroup of $\mathbf{G}$. The correspondence $\rho\mapsto(\rho_I,N_{\rho},\rho_{\Lambda})$ is an order preserving bijection \cite[Theorem 3.5.9]{Howie}. We recall that the normal subgroup $N_{\rho}$ determines the congruence $\rho_G$ on the group $\mathbf G$ in the usual way. Hence, for every $g,h\in G$ we have $g\,\rho_G\,h$ if and only if $(1,g,1) \,\rho\, (1,h,1)$.

\begin{de}(\cite{Howie}, Chapter 2.5)\label{DefnOrt}
A regular semigroup $\bS$ is \emph{orthodox} if its set of idempotents $E(S)$ forms a subsemigroup.
\end{de}

For $k\in\N$, we will denote by $\ox$ a $k$-tuple $(x^1,\dotsc,x^k)$
from $S^k$. Also, if $\alpha$ is a congruence on $\bS$, and
$\ox,\oy\in S^k$, we will use notation $\ox\;\alpha\;\oy$ to denote
that for every $i\in\{1,\dotsc,k\}$ we have $x^i\;\alpha\;y^i$.
Sometimes we put the domain of the semigroup $S$ as the index of the
relation to avoid confusion (for example $1_S$ or $0_S$).

By $\Con(\bA)$ we denote the lattice of all congruences on an algebra $\bA$. The following definition is due to A. Bulatov \cite{Bu:OTNO}, and it presents a generalization of the binary term condition commutator.

\begin{de}(cf.\cite{Bu:OTNO})\label{DefnC(a,a,a,b;d)}
Let $\bA$ be an algebra. Let $k\in\N$ and let $\alpha_1,\dotsc,\alpha_k,\beta$ and $\delta$ be congruences in $\Con(\bA)$. Then we say that $\alpha_1,\dotsc,\alpha_k$ \textit{centralize} $\beta$ \textit{modulo} $\delta$, and we write $C(\alpha_1,\dotsc,\alpha_k,\beta;\delta)$, if for every polynomial $p\in\Pol(\bA)$, of arity $n_1+\dotsc+n_k+m$, where $n_1,\dotsc,n_k,m\in\N$, and every $\oa_1,\ob_1\in A^{n_1}$, $\dotsc$, $\oa_{k},\ob_{k}\in A^{n_{k}}$ and $\oc,\od\in A^{m}$ such that
\begin{enumerate}
\item[(1)] $\oa_j \;\alpha_j\; \ob_j$, for $j=1,\dotsc,k$;
\item[(2)] $\oc \;\beta\; \od$;
\item[(3)] $p(\oa_1,\oa_2,\dotsc,\oa_{k},\oc)\;\delta\;p(\oa_1,\oa_2,\dotsc,\oa_{k},\od)$\\
$p(\oa_1,\ob_2,\dotsc,\oa_{k},\oc)\;\delta\;p(\oa_1,\ob_2,\dotsc,\oa_{k},\od)$\\
\hspace*{30mm}$\dots$\\
$p(\ob_1,\dotsc,\ob_{k-1},\oa_{k},\oc)\;\delta\;p(\ob_1,\dotsc,\ob_{k-1},\oa_{k},\od)$
\end{enumerate}
we have $p(\ob_1,\dotsc,\ob_{k},\oc) \; \delta \; p(\ob_1,\dotsc,\ob_{k},\od)$.
\end{de}

\begin{de}(\cite{Bu:OTNO})\label{DefnHC}
Let $\bA$ be an algebra. Let $k\geq 2$ and let $\alpha_1,\dotsc,\alpha_k\in\Con(\bA)$. The \textit{$k$-ary commutator of $\alpha_1,\dotsc,\alpha_k$} is the smallest congruence $\delta$ on $\bA$ such that $C(\alpha_1,\dotsc,\alpha_k;\delta)$.
\end{de}

The $k$-ary commutator of $\alpha_1,\dotsc,\alpha_k$ will be denoted by $[\alpha_1,\dotsc,\alpha_k]$.
For $k=2$ and congruences $\alpha_1,\alpha_2\in\Con(\bA)$, the commutator $[\alpha_1,\alpha_2]$ defined by Definition \ref{DefnHC} is the binary term condition commutator (see \cite{McMc}, Definition 4.150 and Exercises 4.156.2). This has been proved in \cite{EA:TPFoCA}.

\begin{pro}\textup{(\cite{FM:CTCMV}, Proposition 3.4)} \label{CommutatorProperties}
Let $\bA$ be an algebra, and let $\alpha,\beta,\gamma,\delta$ be congruences on $\bA$. Then we have
\begin{enumerate}
    \item[(i)] $[\alpha,\beta] \leq \alpha\wedge\beta$;
    \item[(ii)] if $\gamma\leq\alpha$ and $\delta\leq\beta$, then $[\gamma,\delta]\leq[\alpha,\beta]$.
\end{enumerate}
\end{pro}

\begin{pro}\label{HC3}
Let $\bA$ be an algebra, let $n\in\N$ and let $\alpha_1,\dots,\alpha_{n}$ be congruences on $\bA$. Then we have $[\alpha_1,\dots,\alpha_n]\leq[\alpha_{n-1},\alpha_n]$.
\end{pro}

\begin{proof} This follows from the property (HC3) in \cite{AM:Magna} and \cite{Bu:OTNO}.
\end{proof}

\begin{de}\textup{(\cite{FM:CTCMV}, Definition 6.1)}\label{DefNilpSolv}
Let $\bA$ be an algebra, and let $\rho,\sigma$ be congruences on $\bA$. We define the series $(\rho,\sigma]^{(k)}$, $k\in\N$, as follows: $(\rho,\sigma]^{(1)} = [\rho,\sigma]$, and $(\rho,\sigma]^{(k+1)} = [\rho,(\rho,\sigma]^{(k)}]$, for $k\in\N$. Similarly, we define the series $[\rho]^{(k)}, k\in\N$ by $[\rho]^{(1)}=[\rho,\rho]$ and $[\rho]^{(k+1)}=[[\rho]^{(k)},[\rho]^{(k)}]$, for $k\in\N$.

Algebra $\bA$ is \textit{$n$-nilpotent}, $n\in\N$ if $(1_A,1_A]^{(n)}=0_A$. Similarly, an algebra $\bA$ is $n$-solvable, $n\in\N$, if $[1_A]^{(n)}=0_A$.
\end{de}

\begin{de}\label{supernilpotent}
Let $k\in\N$. An algebra $\bA$ is $k$-\emph{supernilpotent} if $[\,\underbrace{1_A,\dotsc,1_A}_k\,]=0_A$.
\end{de}

\begin{rem}\label{remNilp}
An $n$-nilpotent algebra is often called \textit{left} $n$-nilpotent. \emph{Right nilpotent} algebras are defined using a dually defined series of congruences. In groups, the notions of left and right nilpotency do coincide. However, for an arbitrary algebra the degrees of left and right nilpotency do not have to be equal (see \cite{KK}). In this paper we will study only the notion of left nilpotency. However, it should be noted that all results obtained for (left) nilpotent semigroups can be dually proved for right nilpotent semigroups.  \\
Also note that the properties of the commutator stated in Proposition \ref{CommutatorProperties} imply the inequality $[\rho]^{(k+1)}\leq (\rho,\rho]^{(k)}$, for every $k\in\N$, and for arbitrary congruence $\rho\in\Con(\bA)$. Hence for $n\in\N$, if an algebra $\bA$ is $n$-nilpotent (left or right), then it is $n$-solvable.
\end{rem}

Let us recall that if we have $n\in\N$, algebras $\mathbf{A}_1,\dots,\mathbf{A}_n$ of the same type and $\theta_1,\dots,\theta_n$ their congruences, respectively, then the direct product of the congruences $\theta_1,\dotsc,\theta_n$, in abbreviation $\theta_1\times\dots\times\theta_n$, is a congruence on $\mathbf{A}_1\times\dots\times\mathbf{A}_n$ given by
$$
(a_1,\dots,a_n)\,\theta_1\times\dots\times\theta_n\,(b_1,\dots,b_n) \Leftrightarrow a_1\,\theta_1\,b_1\wedge\dots\wedge a_n\,\theta_n\,b_n,
$$
for all $ (a_1,\dots,a_n),(b_1,\dots,b_n)\in A_1\times\dots\times A_n$. As usual we will call an algebra $\mathbf{A}_1\times\dots\times\mathbf{A}_n$ skew-free if all its congruences are direct product of congruences of algebras $\mathbf{A}_1,\dots,\mathbf{A}_n$, respectively.

\begin{pro}\label{congproduct}
Every congruence $\rho$ on a completely simple semigroup $\bS=\mathcal{M}[G;I,\Lambda;P]$ satisfies $\rho=\rho_I\times\rho_G\times\rho_{\Lambda}$.
\end{pro}

\begin{proof}
Let $\rho$ be a congruence on the semigroup $\bS$, and let $(i,g,\lambda),(j,h,\mu)\in S$. Assume that $(i,g,\lambda)\,\rho\,(j,h,\mu)$. Then \cite[Lemma 3.5.3, 3.5.4]{Howie} implies $i\,\rho_I\,j$ and $\lambda\,\rho_{\Lambda}\,\mu$. If we multiply the relation $(i,g,\lambda)\;\rho\;(j,h,\mu)$ with $(1,e,1)$ on the left and on the right, we obtain $(1,p_{1i}gp_{\lambda 1},1)\;\rho\;(1,p_{1j}hp_{\mu 1},1)$. By definition of the congruence $\rho_G$, it follows that $p_{1i}gp_{\lambda 1}\;\rho_G\;p_{1j}hp_{\mu 1}$. However, since $p_{1x}=e=p_{\xi 1}$ for every $x\in I$, $\xi\in \Lambda$, it follows that $g\;\rho_G\;h$. Therefore, $(i,g,\lambda)\, \rho_I\times\rho_G\times\rho_{\Lambda}\, (j,h,\mu)$. Now assume that $i\;\rho_I\;j$, $\lambda\;\rho_{\Lambda}\;\mu$ and $g\;\rho_G\;h$. From definition of $\rho_I$ and $\rho_{\Lambda}$, it follows that $(i,e,1)\;\rho\;(j,e,1)$ and $(1,e,\lambda)\;\rho\;(1,e,\mu)$. Since $g\;\rho_G\;h$ implies $(1,g,1)\;\rho\;(1,h,1)$, using $p_{11}=e$ and the compatibility of the congruence $\rho$, we obtain $(i,g,\lambda) = (i,e,1)(1,g,1)(1,e,\lambda) \;\rho\; (j,e,1)(1,h,1)(1,e,\mu) = (j,h,\mu)$.
\end{proof}

\begin{pro}\label{comutatorofproducts}
If $\mathbf{A}$, $\mathbf{B}$ and $\mathbf{C}$ are three algebras of
the same type such that $\bA\times\mathbf{B}\times\mathbf{C}$ is
skew-free, $n\in\N$ and $\alpha_1,\dots,\alpha_n\in\Con{\bA}$,
$\beta_1,\dots,\beta_n\in\Con{\mathbf B}$ and
$\gamma_1,\dots,\gamma_n\in\Con{\mathbf C}$ then,
$$
[\alpha_1\times\beta_1\times\gamma_1,\dots,\alpha_n\times\beta_n\times\gamma_n]=[\alpha_1,\dots,\alpha_n]\times[\beta_1,\dots,\beta_n]\times[\gamma_1,\dots,\gamma_n].
$$
\end{pro}

\begin{proof} By Definition \ref{DefnHC} and definition of skew-free algebras and direct product of congruences.
\end{proof}

\begin{pro}\textup{(\cite[Theorem 6.8]{AE:IJAC}, \cite[Corollary 6.15]{AM:Magna})}\label{GroupNilpSupernilp}
Let $\bG$ be a group, and let $n\in\N$. Then $\bG$ is $n$-supernilpotent if and only if $\bG$ is $n$-nilpotent.
\end{pro}

\section{Regular semigroups}

\begin{pro}\textup{(\cite{Warne}, Corollary 2.6)}\label{WarneCor2.6}
A regular semigroup $\bS$ is abelian if and only if $\bS$ is isomorphic to the direct product $\mathbf{G}\times \mathbf{L}\times \mathbf{R}$, where $\mathbf{G}$ is an abelian group, $\mathbf{L}$ is a left zero semigroup and $\mathbf{R}$ is a right zero semigroup.
\end{pro}

\begin{pro}\label{ProIdempNilpSolv}
Let $n\geq 2$, let $\bS$ be a semigroup with congruences $\alpha_1,\dots,\alpha_n$ and let $e$ and $f$ be idempotents in $\bS$ such that $e\,\alpha_i\, f$ for all $i\in\{1,\dots,n\}$. If $e\leq f$, then we have $e\,[\alpha_1,\dots,\alpha_n]\, f$.
\end{pro}

\begin{proof}
We take the polynomial $p(x_1,\dots,x_n)=x_1\cdot\ldots\cdot x_n$. Inequality $e\leq f$ gives us $ef=e=fe$. Hence,
$p(x_1,\dots,x_{n-1},e)=e=p(x_1,\dots,x_{n-1},f)$ for all $(x_1,\dots,x_{n-1})\in\{e,f\}^{n-1}\backslash\{(f,\dots,f)\}$. Therefore, $p(x_1,\dots,x_{n-1},e)\,[\alpha_1,\dots,\alpha_n]$ $p(x_1,\dots,x_{n-1},f)$ for all $(x_1,\dots,x_{n-1})\in\{e,f\}^{n-1}\backslash\{(f,\dots,f)\}$. Since we have $C(\alpha_1,\dots,\alpha_n;[\alpha_1,\dots,\alpha_n])$ we obtain $e=p(f,\dots,f,e)\,[\alpha_1,\dots,\alpha_n]\,p(f,\dots,f,f)=f$.
\end{proof}

\begin{cor}\label{Cor:e<=f}
Let $\bS$ be a semigroup. If the semigroup $\bS$ satisfies either of the following conditions
\begin{enumerate}
\item[(i)] $\bS$ is $n$-nilpotent,
\item[(ii)] $\bS$ is $n$-solvable,
\item[(iii)] $\bS$ is $n$-supernilpotent,
\end{enumerate}
for some $n\in\N$, then for every $e,f\in E(S)$ the inequality $e\leq f$ implies $e=f$, that is, $(E(S),\leq)$ is an antichain.
\end{cor}

\begin{proof}
Since $e,f\in S$ are idempotents such that $e\leq f$ each of the assumptions (i)-(iii) implies the statement by Proposition \ref{ProIdempNilpSolv} and definition of nilpotency, solvability and supernilpotency, because $e\,1\,f$.
\end{proof}

\begin{cor}\cite[Theorem 3.3.3 $(4)\Rightarrow(1)$]{Howie}\label{CorOfLemmaID}
Let $\bS$ be a regular semigroup such that $(E(S),\leq)$ is an antichain, where $\leq$ is the natural partial ordering of idempotents. Then $\bS$ is a completely simple semigroup.
\end{cor}

\begin{pro}\label{ThmRegNilpSolv}
Let $\bS$ be a regular semigroup, and let $n\in\N$.
\begin{enumerate}
\item[(i)] The semigroup $\bS$ is $n$-nilpotent if and only if $\bS$ is an $n$-nilpotent completely simple semigroup.
\item[(ii)] The semigroup $\bS$ is $n$-solvable if and only if $\bS$ is an $n$-solvable completely simple semigroup.
\item[(iii)] The semigroup $\bS$ is $n$-supernilpotent if and only if $\bS$ is an $n$-supernilpotent completely simple semigroup.
\end{enumerate}
\end{pro}

\begin{proof}
(i), (ii) and (iii) ($\leftarrow$) Every completely simple semigroup is regular, hence this implication is trivially true in statements (i), (ii) and (iii).

\noindent (i), (ii) and (iii) ($\rightarrow$) Assume that $\bS$ has one of the property mentioned in (i), (ii) or (iii). From Corollary \ref{Cor:e<=f}, it follows that $(E(S),\leq)$ is an antichain. Therefore, by Corollary \ref{CorOfLemmaID} we obtain that $\bS$ is a completely simple semigroup with the same property.
\end{proof}

Let us notice here that neither nilpotency nor supernilpotency in regular semigroups does not carry out for the semigroup to be orthodox. We provide the following example.

\begin{exa}\label{nilpotentnotorthodox}
We take the cyclic group of order $C_2=\{e,g\}$ and Rees matrix semigroup $S_2=\mathcal{M}[C_2;\{1,2\},\{1,2\};P]$ where $P=\left[\begin{array}{cc}e&e\\e&g\end{array}\right]$. $S_2$ is  regular because it is completely simple semigroup. One can easily check that $(1,e,1)$ and $(2,g,2)$ are idempotents, but
$$
(1,e,1)(2,g,2)=(1,g,2)\neq(1,e,2)=((1,e,1)(2,g,2))^2.
$$
Hence, $S_2$ is not an orthodox semigroup. Therefore, $S_2$ is not abelian by Proposition \ref{WarneCor2.6}. In Propositions \ref{2nilpotent} and \ref{2supernilpotent} we prove that it is $2$-nilpotent and $2$-supernilpotent.
\end{exa}



Note that due to the commutativity of the group $C_2$, for any polynomial $t\in
\Pol(\bS)$, $\ar(t)=n_1+\dotsc+n_k$, where $n_1,\dotsc,n_k\in\N$,
and any $\ox_s\in S_2^{n_s}$, for $s=1,\dotsc,k$, there exist
$\ell\in I$, $\alpha\in\N_0$ and $\sigma\in\Lambda$ such that
\begin{equation*}
t(\ox_1,\dotsc,\ox_k) = (\ell,
g^{\alpha}\cdot t^{C_2}(\ox_1,\dotsc,\ox_k),\sigma)
\end{equation*}
where $\alpha$ is the number of all $p_{22}=g$ which appear in the
group coordinate of the product $t(\ox_1,\dotsc,\ox_k)$. In the
following proofs we will use this property without explicit
referencing.

\begin{lem}\label{LemaPom}
Let $\bS=\mathcal{M}[G;I,\Lambda;P]$ be a completely simple semigroup, and let $\delta\in\Con(\bS)$ be such that $\delta_I=0_I$ and $\delta_{\Lambda}=0_{\Lambda}$. Let $t$ be a polynomial over $\bS$, $\ar(t)=n_1+\dotsc+n_{k-1}+m$, where $n_1,\dotsc,n_{k-1},m\in\N$, and let $\oa_s,\ob_s\in S_2^{n_s}$, for $s=1,\dotsc,k-1$, $\oc,\od\in S_2^m$ be such that
\begin{align}
t(\oa_1,\oa_2,\dotsc,\oa_{k-1},\oc) &\;\delta\; t(\oa_1,\oa_2,\dotsc,\oa_{k-1},\od) \nonumber \\
t(\oa_1,\ob_2,\dotsc,\oa_{k-1},\oc) &\;\delta\; t(\oa_1,\ob_2,\dotsc,\oa_{k-1},\od) \nonumber \\
&\dots \nonumber \\
t(\ob_1,\dotsc,\ob_{k-2},\oa_{k-1},\oc) &\;\delta\; t(\ob_1,\dotsc,\ob_{k-2},\oa_{k-1},\od)\label{eqn:LemaPom}
\end{align}
If $(\ell,h,\sigma) = t(\ob_1,\dotsc,\ob_{k-1},\oc)$ and $(\ell',h',\sigma')=t(\ob_1,\dotsc,\ob_{k-1},\od)$, then $\ell=\ell'$ and $\sigma=\sigma'$.
\end{lem}

\begin{proof}
We will prove that $\ell=\ell'$, while $\sigma=\sigma'$ is proved
dually. We will assume that the set of all $2^{k-1}$ $(k-1)$-tuples
$(\ox_1,\dotsc,\ox_{k-1})$, where $\ox_s\in\{\oa_s,\ob_s\}$,
$s=1,\dotsc,k-1$ is ordered using the lexicographic order. Here we
use the notation $a^u_s=(i^u_s,f^u_s,\lambda^u_s)$,
$b^u_s=(i'^u_s,f'^u_s,\lambda'^u_s)$, for $u=1,\dotsc,n_s$, while
$c^v=(j_v,g_v,\mu_v)$, $d^v=(k_v,h_v,\nu_v)$, $v=1,\dotsc,m$ and
introduce the following:
\begin{align*}
&(\ell_s,h_s,\sigma_s) = t(\ox_1,\dotsc,\ox_{k-1},\oc)\\
&(\ell_s',h_s',\sigma_s') = t(\ox_1,\dotsc,\ox_{k-1},\od)
\end{align*}
where $(\ox_1,\dotsc,\ox_{k-1})$ is the $s$-th $(k-1)$-tuple in the ordering. Let us notice that $\ell=\ell_{2^{k-1}}$ and $\ell'=\ell_{2^{k-1}}'$. We differentiate between three possibilities:\\
\textbf{case 1:} $t(\ox_1,\dotsc,\ox_{k-1},\oy) = c \cdot t_1(\ox_1,\dotsc,\ox_{k-1},\oy)$ where $c=(j,g,\mu)$ is a constant from $\bS$. Then we have $\ell_s=j=\ell_s'$ for every $s$, and consequently $\ell=j=\ell'$. \\
\textbf{case 2:} $t(\ox_1,\dotsc,\ox_{k-1},\oy) = x_p^r\cdot t_1(\ox_1,\dotsc,\ox_{k-1},\oy)$ for some $p\in\{1,\dotsc,k-1\}$ and $r\in\{1,\dotsc,n_p\}$, and some $t_1\in\Pol(\bS)$. Then in both $t(\ob_1,\dotsc,\ob_{k-1},\oc)$ and $t(\ob_1,\dotsc,\ob_{k-1},\od)$ we have $\ob_p$ on the $p$-th place with coordinate $b^r_p=(i'^r_p,f'^r_p,\lambda'^r_p)$ and therefore $\ell=i'^r_p=\ell'$.\\
\textbf{case 3:} $t(x,\oy) = y_i\cdot
t_1(\ox_1,\dotsc,\ox_{k-1},\oy)$ for some $i\in\{1,\dotsc,m\}$ and
some $t_1\in\Pol(\bS)$. In this case we have $\ell_s=\mu_i$ and
$\ell_s'=\nu_i$ for every $s$. Since $\delta_I=0_I$, relations
\eqref{eqn:LemaPom} imply that $\mu_i=\nu_i$, therefore
$\ell=\mu_i=\nu_i=\ell'$, as claimed.
\end{proof}

\begin{lem}\label{1,1}
$C(1_{S_2},1_{S_2};\rho)$.
\end{lem}

\begin{proof}
Let $t\in\Pol_{n+m}(S_2)$, $n,m\in\N$ and let $\oa,\ob\in S_2^n$,
$\oc,\od\in S_2^m$ be such that $t(\oa,\oc)\;\rho\;t(\oa,\od)$. We
will introduce the following notation
$t(\oa,\oc)=(\ell_1,h_1,\sigma_1)$,
$t(\oa,\od)=(\ell_1',h_1',\sigma_1')$, $t(\ob,\oc)=(\ell_2,h_2,\sigma_2)$ and
$t(\ob,\od)=(\ell_2',h_2',\sigma_2')$. Since $\rho_I=0_I$ and
$\rho_{\Lambda}=0_{\Lambda}$, by Lemma \ref{LemaPom} from relation
$t(\oa,\oc)\;\rho\;t(\oa,\od)$ it follows that $\ell_2=\ell_2'$ and
$\sigma_2=\sigma'_2$. On the other hand, since $N_{\rho}=C_2=G$, it
follows that $h_2\;\rho_G\;h_2'$ is trivially true.  Therefore, we
have $t(\ob,\oc)\;\rho\;t(\ob,\od)$, by Proposition \ref{congproduct}. This proves the centralizing
condition $C(1_{S_2},1_{S_2};\rho)$.
\end{proof}

\begin{lem}\label{1,rho}
$C(1_{S_2},\rho;0_{S_2})$.
\end{lem}

\begin{proof}
Let $t\in\Pol_{n+m}(S_2)$, $n,m\in\N$ and let $\oa,\ob\in S_2^n$, $\oc,\od\in S_2^m$ be such that $\oc\;\rho\;\od$.
Here we use the notation $a^u=(i_u,f_u,\lambda_u)$, $b^u=(i_u',f_u',\lambda_u')$, for $u=1,\dotsc,n$, while $c^v=(j_v,g_v,\mu_v)$, $d^v=(k_v,h_v,\nu_v)$, $v=1,\dotsc,m$. Also denote $\mathbf{f}=(f_1,\dotsc,f_n)$ and $\mathbf{f'}=(f_1',\dotsc,f_n')$, as well as $\mathbf{g}=(g_1,\dotsc,g_m)$ and $\mathbf{h}=(h_1,\dotsc,h_m)$. From $c^v\;\rho\;d^v$ it follows that $j_v=k_v$, $\mu_v=\nu_v$ and $g_vh_v^{-1}\in C_2$. The last condition is trivially satisfied, and it can be written as $h_v=g^{\epsilon_v}g_v$ where $\epsilon_v\in\{0,1\}$ for $v=1,\dotsc,m$. We will also use the notation
\begin{align*}
& t(\oa,\oc) = (\ell_1,g^{\alpha_1}\cdot t^{C_2}(\mathbf{f},\mathbf{g}),\sigma_1) \mbox{ and } t(\oa,\od) = (\ell_2,g^{\alpha_2} \cdot t^{C_2}(\mathbf{f},\mathbf{h}),\sigma_2);\\
& t(\ob,\oc) = (\ell_1',g^{\alpha_1'} \cdot t^{C_2}(\mathbf{f}',\mathbf{g}),\sigma_1') \mbox{ and } t(\ob,\od) = (\ell_2',g^{\alpha_2'} \cdot t^{C_2}(\mathbf{f}',\mathbf{h}),\sigma_2'),
\end{align*}
for some $\ell_s,\ell_s'\in I$, $\sigma_s,\sigma_s'\in\Lambda$, $s=1,2$, where $\alpha_s,\alpha_s'$, $s=1,2$ denote the number of $p_{\lambda i}$-s equal to $g$ in the group part of the corresponding product. Since $j_v=k_v$ and $\mu_v=\nu_v$, for every $v=1,\dotsc,m$, it follows that $\alpha_1=\alpha_2$ and $\alpha_1'=\alpha_2'$.

Assume that $t(\oa,\oc)=t(\oa,\od)$, then we have $\ell_1=\ell_2$, $\sigma_1=\sigma_2$ and $g^{\alpha_1}\cdot t^{C_2}(\mathbf{f},\mathbf{g})=g^{\alpha_2}\cdot t^{C_2}(\mathbf{f},\mathbf{h})$. Since $\alpha_1=\alpha_2$, it follows that $t(\mathbf{f},\mathbf{g})=t(\mathbf{f},\mathbf{h})$.  Using the equalities $h_v=g^{\epsilon_v}g_v$, $v=1,\dotsc,m$ and the commutativity of $C_2$, we obtain
\begin{align*}
&t(\mathbf{f},\mathbf{h}) = t(e,\dotsc,e,g^{\epsilon_1},\dotsc,g^{\epsilon_n})\cdot t(\mathbf{f},\mathbf{g}),\\ &t(\mathbf{f'},\mathbf{h}) = t(e,\dotsc,e,g^{\epsilon_1},\dotsc,g^{\epsilon_n})\cdot t(f',\mathbf{g}).
\end{align*}
Now the equality $t(\mathbf{f},\mathbf{g})=t(\mathbf{f},\mathbf{h})$ implies $e=t(e,\dotsc,e,g^{\epsilon_1},\dotsc,g^{\epsilon_n})$. Hence, we have $t(f',\mathbf{h}) = t(\mathbf{f}',\mathbf{g})$. Since we also have $\alpha_1'=\alpha_2'$, the previous equality gives us $g^{\alpha_1'}\cdot t(\mathbf{f'},\mathbf{g}) = g^{\alpha_2'}\cdot t(\mathbf{f'},\mathbf{h})$. By Lemma \ref{LemaPom} we also obtain $\ell_1'=\ell_2'$ and $\sigma_1'=\sigma_2'$, because $=_I=0_I$ and $=_{\Lambda}=0_{\Lambda}$. Hence, we have $t(\ob,\oc)=t(\ob,\od)$, as claimed.
\end{proof}

\begin{pro}\label{2nilpotent}
$S_2$ is $2$-nilpotent.
\end{pro}

\begin{proof}
We obtain $[1_{S_2},[1_{S_2},1_{S_2}]]\leq[1_{S_2},\rho]=0_{S_2}$, using Proposition \ref{CommutatorProperties}, Lemma \ref{1,1} and Lemma \ref{1,rho}, and therefore 
$S_2$ is 2-nilpotent.
\end{proof}

\begin{pro}\label{2supernilpotent}
$S_2$ is $2$-supernilpotent.
\end{pro}

\begin{proof}
According to the Definition \ref{DefnHC}, it is enough to prove the centralizing condition
$C(1_{S_2},1_{S_2},1_{S_2};0_{S_2})$. Let $t\in\Pol(S_2)$, $\ar(t)=n_1+n_2+m$, and let $\oa_1,\ob_1\in S_2^{n_1}$, $\oa_2,\ob_2\in S_2^{n_2}$, $\oc,\od\in S_2^m$. Assume that
\begin{align}
& t(\oa_1,\oa_2,\oc)=t(\oa_1,\oa_2,\od) \nonumber \\
& t(\oa_1,\ob_2,\oc)=t(\oa_1,\ob_2,\od) \nonumber \\
& t(\ob_1,\oa_2,\oc)=t(\ob_1,\oa_2,\od)\label{eqnTC}
\end{align}
Here we are using the following notation $a_1^u=(i_u,f_u,\lambda_u)$, $b_1^u=(i_u',f_u',\lambda_u')$ for $u=1,\dotsc,n_1$, $a_2^u=(i_u,f_u,\lambda_u)$, $b_2^u=(i_u',f_u',\lambda_u')$ for $u=n_1,\dotsc,n_1+n_2$, and $c^v=(j_v,g_v,\mu_v)$, $d^v=(k_v,h_v,\nu_v)$ for $v=1,\dotsc,m$. We will also denote $\mathbf{f_1}=(f_1,\dots,f_{n_1})$, $\mathbf{f_1}'=(f_1',\dotsc,f_{n_1}')$, and $\mathbf{f_2}=(f_{n_1},\dots,f_{n_1+n_2})$, $\mathbf{f_2}'=(f_{n_1}',\dots,f_{n_1+n_2}')$, as well as $\mathbf{g}=(g_1,\dotsc,g_m)$ and $\mathbf{h}=(h_1,\dots,h_m)$. Since the group $C_2$ is commutative, there exist polynomials $p,q,r\in\Pol(C_2)$ such that
\begin{equation*}
t^{C_2}(\ox,\oy,\oz) = p(\ox)\cdot q(\oy) \cdot r(\oz).
\end{equation*}
Therefore, we can write
\begin{align*}
& t(\oa_1,\oa_2,\oc) = (\ell_1,g^{\alpha_1}\cdot t^{C_2}(\mathbf{f}_1,\mathbf{f}_2,\mathbf{g}),\sigma_1) = (\ell_1,g^{\alpha_1} \cdot p(\mathbf{f}_1)\cdot q(\mathbf{f}_2)\cdot r(\mathbf{g}),\sigma_1) \\
& t(\oa_1,\oa_2,\od) = (\ell_2,g^{\alpha_2}\cdot t^{C_2}(\mathbf{f}_1,\mathbf{f}_2,\mathbf{h}),\sigma_2) = (\ell_2,g^{\alpha_2} \cdot p(\mathbf{f}_1)\cdot q(\mathbf{f}_2)\cdot r(\mathbf{h}),\sigma_2)\\
& t(\oa_1,\ob_2,\oc) = (\ell_3,g^{\alpha_3}\cdot t^{C_2}(\mathbf{f}_1,\mathbf{f}_2',\mathbf{g}),\sigma_3) = (\ell_3,g^{\alpha_3} \cdot p(\mathbf{f}_1) \cdot q(\mathbf{f}_2') \cdot r(\mathbf{g}),\sigma_3) \\
& t(\oa_1,\ob_2,\od) = (\ell_4,g^{\alpha_4}\cdot t^{C_2}(\mathbf{f}_1,\mathbf{f}_2',\mathbf{h}),\sigma_4) = (\ell_4,g^{\alpha_4} \cdot p(\mathbf{f}_1) \cdot q(\mathbf{f}_2') \cdot r(\mathbf{h}),\sigma_4)\\
& t(\ob_1,\oa_2,\oc) = (\ell_5,g^{\alpha_5}\cdot t^{C_2}(\mathbf{f}_1',\mathbf{f}_2,\mathbf{g}),\sigma_5) = (\ell_5,g^{\alpha_5} \cdot p(\mathbf{f}_1')\cdot q(\mathbf{f}_2) \cdot r(\mathbf{g}),\sigma_5)\\
& t(\ob_1,\oa_2,\od) = (\ell_6,g^{\alpha_6}\cdot t^{C_2}(\mathbf{f}_1',\mathbf{f}_2,\mathbf{h}),\sigma_6) = (\ell_6,g^{\alpha_6} \cdot p(\mathbf{f}_1') \cdot q(\mathbf{f}_2) \cdot r(\mathbf{h}),\sigma_6)\\
& t(\ob_1,\ob_2,\oc) = (\ell_7,g^{\alpha_7}\cdot t^{C_2}(\mathbf{f}_1',\mathbf{f}_2',\mathbf{g}),\sigma_7) = (\ell_7,g^{\alpha_7}\cdot p(\mathbf{f}_1') \cdot q(\mathbf{f}_2') \cdot r(\mathbf{g}),\sigma_7)\\
& t(\ob_1,\ob_2,\od)=(\ell_8,g^{\alpha_8}\cdot t^{C_2}(\mathbf{f}_1',\mathbf{f}_2',\mathbf{h}),\sigma_8) = (\ell_8,g^{\alpha_8} \cdot p(\mathbf{f}_1') \cdot q(\mathbf{f}_2')
\cdot r(\mathbf{h}),\sigma_8)
\end{align*}
Now from the equalities \eqref{eqnTC} we have
\begin{align*}
& g^{\alpha_1}\cdot  p(\mathbf{f}_1)\cdot q(\mathbf{f}_2) \cdot r(\mathbf{g}) = g^{\alpha_2} \cdot p(\mathbf{f}_1) \cdot q(\mathbf{f}_2) \cdot r(\mathbf{h})\\
& g^{\alpha_3} \cdot p(\mathbf{f}_1) \cdot q(\mathbf{f}_2')\cdot r(\mathbf{g}) = g^{\alpha_4} \cdot p(\mathbf{f}_1) \cdot q(\mathbf{f}_2') \cdot r(\mathbf{h})\\
& g^{\alpha_5} \cdot p(\mathbf{f}_1')\cdot q(\mathbf{f}_2) \cdot r(\mathbf{g}) = g^{\alpha_6} \cdot  p(\mathbf{f}_1') \cdot q(\mathbf{f}_2)\cdot r(\mathbf{h})
\end{align*}
Multiplying these three equalities we obtain
\begin{align*}
g^{\alpha_1+\alpha_3+\alpha_5} \cdot (p(\mathbf{f}_1))^2 (q(\mathbf{f}_2))^2 p(\mathbf{f}_1') q(\mathbf{f}_2') (r(\mathbf{g}))^3  \\
=  g^{\alpha_2+\alpha_4+\alpha_6} \cdot (p(\mathbf{f}_1))^2 (q(\mathbf{f}_2))^2 p(\mathbf{f}_1') q(\mathbf{f}_2') (r(\mathbf{h}))^3.
\end{align*}
Hence $g^{\alpha_1+\alpha_3+\alpha_5} \cdot p(\mathbf{f}_1') q(\mathbf{f}_2') r(\mathbf{g}) = g^{\alpha_2+\alpha_4+\alpha_6} \cdot p(\mathbf{f}_1') q(\mathbf{f}_2') r(\mathbf{h})$, because $C_2$ is a group of order $2$. This further implies
\begin{equation}
g^{\alpha_1+\alpha_3+\alpha_5 - \alpha_2 -\alpha_4-\alpha_6} \cdot t^{C_2}(\mathbf{f}_1',\mathbf{f}_2',\mathbf{g}) = t^{C_2}(\mathbf{f}_1',\mathbf{f}_2',\mathbf{h}).\label{eqn:Eq135}
\end{equation}
We will denote by $\epsilon(\ox,\oy)$ the number of $p_{22}=g$-s that appear in the product of a coordinate from $\ox$ with a coordinate from $\oy$, in any order. In particular, $\epsilon(\ox,const)$ is the number of $p_{22}=g$-s that appear in the product of a coordinate of $\ox$ and a constant, also in any order. Using this notation, we obtain
\begin{align*}
           & \alpha_1+\alpha_3+\alpha_5-\alpha_2-\alpha_4-\alpha_6 = (\alpha_1-\alpha_2) + (\alpha_3-\alpha_4) + (\alpha_5-\alpha_6) \\
\equiv \;  & (\epsilon(\oa_1,\oc)+\epsilon(\oa_2,\oc)+\epsilon(\oc,const)+\epsilon(\oc,\oc)) \\
          -& (\epsilon(\oa_1,\od)+\epsilon(\oa_2,\od)+\epsilon(\od,const)+\epsilon(\od,\od))\\
          +& (\epsilon(\oa_1,\oc)+\epsilon(\ob_2,\oc)+\epsilon(\oc,const)+\epsilon(\oc,\oc))\\
          -& (\epsilon(\oa_1,\od)+\epsilon(\ob_2,\od)+\epsilon(\od,const)+\epsilon(\od,\od))\\
          +& (\epsilon(\ob_1,\oc)+\epsilon(\oa_2,\oc)+\epsilon(\oc,const)+\epsilon(\oc,\oc))\\
          -& (\epsilon(\ob_1,\od)+\epsilon(\oa_2,\od)+\epsilon(\od,const)+\epsilon(\od,\od))\\
\equiv \;  & 2(\epsilon(\oa_1,\oc)-\epsilon(\oa_1,\od)) + 2(\epsilon(\oa_2,\oc)-\epsilon(\oa_2,\od)) \\
          +& (\epsilon(\ob_1,\oc) - \epsilon(\ob_1,\od)) + (\epsilon(\ob_2,\oc)-\epsilon(\ob_2,\od)) \\
          +& 3(\epsilon(\oc,const)-\epsilon(\od,const)) + 3(\epsilon(\oc,\oc)-\epsilon(\od,\od)) \\
\equiv \;  & (\epsilon(\ob_1,\oc) + \epsilon(\ob_2,\oc) + \epsilon(\oc,const) + \epsilon(\oc,\oc)) \\
          -& (\epsilon(\ob_1,\od)+\epsilon(\ob_2,\oc)+\epsilon(\od,const)+\epsilon(\od,\od))\\
\equiv \;  & \alpha_7-\alpha_8  \pmod 2 
\end{align*}
Since $g$ is an element of order $2$, it follows that $g^{\alpha_1+\alpha_3+\alpha_5-\alpha_2-\alpha_4-\alpha_6} = g^{\alpha_7-\alpha_8}$. Therefore, from the equality \eqref{eqn:Eq135} we obtain
\begin{equation*}
g^{\alpha_7} \cdot t^{C_2}(\mathbf{f}_1',\mathbf{f}_2',\mathbf{g}) = g^{\alpha_8} \cdot t^{C_2}(\mathbf{f}_1',\mathbf{f}_2',\mathbf{h}).
\end{equation*}
Further on, by Lemma \ref{LemaPom}, the equalities \eqref{eqnTC} imply $\ell_7=\ell_8$ and $\sigma_7=\sigma_8$, because $=_I=0_I$ and $=_{\Lambda}=0_{\Lambda}$. Hence, we have proved that
\begin{equation*}
t(\ob_1,\ob_2,\oc) = (\ell_7,g^{\alpha_7} \cdot t^{C_2}(\mathbf{f}_1',\mathbf{f}_2',\mathbf{g}), \sigma_7) = (\ell_8,g^{\alpha_8} \cdot t^{C_2}(\mathbf{f}_1',\mathbf{f}_2',\mathbf{h}), \sigma_8) = t(\ob_1,\ob_2,\od).
\end{equation*}
which completes the proof of $C(1_{S_2},1_{S_2},1_{S_2};0_{S_2})$.
\end{proof}

\section{Orthodox semigroups}

One characterization of completely simple orthodox semigroups is given as an exercise in \cite[Exercise 4.10]{Howie}, and its full proof can be found in \cite[Exercise 5.6]{ACain}. It does not give us the descriptions of the isomorphic structures in detail. Hence we will be using modified versions of that result, which can again be found as exercises in \cite{Howie} and \cite{ACain}.

\begin{pro}\textup{(\cite{Howie}, Exercise 3.8)}\label{HowieExer3.8}
Let $\bS=\mathcal{M}[G;I,\Lambda;P]$ be a completely simple orthodox semigroup. Then $\bS$ is isomorphic to the direct product $\bG\times\mathbf{B}$, where $\mathbf{B}=I\times \Lambda$ is a rectangular band.
\end{pro}

\begin{pro} \textup{(cf. \cite[Exercise 4.10]{Howie}, \cite{ACain}, Exercise 5.6)}\label{PropBandTimesGroup}
Let $\bG$ be a group, let $\mathbf{B}=I\times\Lambda$ be a rectangular band, and let the semigroup $\bS$ be the direct product $\bG\times\mathbf{B}$. Then $\bS$ is an orthodox completely simple semigroup, isomorphic to the Rees matrix semigroup $\mathcal{M}[G;I,\Lambda;P]$, where $P=[e]_{\Lambda\times I}$ is the $\Lambda\times I$ matrix whose all the entries are the identity $e$ of the group $\bG$.
\end{pro}

The characterization of abelian regular semigroups obtained by Warne can be used to deduce the following result easily.

\begin{cor}\label{CorollaryofWarne}
All left(right)-zero semigroups are abelian.
\end{cor}

\begin{proof}
Let $\bS$ be a left(right)-zero semigroup. Then we know that it is isomorphic to the direct product $\mathbf{S}\times\{e\}\times\{e\}$, where $\{e\}$ is trivially an abelian group and also right-zero semigroup. Therefore, $\bS$ is abelian by Proposition \ref{WarneCor2.6}.
\end{proof}

\begin{pro}\label{GIL}
Let $\bS=\bG\times \mathbf{I}\times\mathbf{\Lambda}$, where $\bG$ is a group, $\mathbf{I}$ is a left-zero semigroup and $\mathbf{\Lambda}$ is a right-zero semigroup. Then for every $k\in\N$, $\bS$ is
\begin{enumerate}
\item $k$-nilpotent;
\item $k$-solvable;
\item $k$-supernilpotent
\end{enumerate}
 if and only if the same property is true for $\bG$.
\end{pro}

\begin{proof}
Let $k\in\N$. We prove the item (3) and the other two items are analogous. Since $\bS$ is completely simple semigroup by Proposition \ref{PropBandTimesGroup}, it is also skew-free by Proposition \ref{congproduct}. Hence, $\bS$ is $k$-supernilpotent if and only if $\mathbf I$, $\mathbf\Lambda$ and $\bG$ are $k$-supernilpotent, by Proposition \ref{comutatorofproducts}. By Corollary \ref{CorollaryofWarne} we know that left and right-zero semigroups are always abelian and hence $k$-supernilpotent using Proposition \ref{HC3}. Therefore, we obtain that $\bS$ is $k$-supernilpotent if and only if $\bG$ is $k$-supernilpotent.
\end{proof}

\begin{cor}\label{SuperNilpOrthCSsg}
Let $\bS=\bG\times \mathbf{I}\times\mathbf{\Lambda}$, where $\bG$ is a group, $\mathbf{I}$ is a left-zero semigroup and $\mathbf{\Lambda}$ is a right-zero semigroup. Then for every $k\in\N$, $\bS$ is $k$-supernilpotent if and only if $\bG$ is $k$-nilpotent.
\end{cor}

\begin{proof}
Using Propositions \ref{GIL} and \ref{GroupNilpSupernilp} we obtain the statement.
\end{proof}

\begin{thm}\label{SupernilpOrthodox}
Let $\bS$ be an orthodox semigroup, and let $n\in\N$. Then $\bS$ is $n$-supernilpotent if and only if $\bS$ is isomorphic to a direct product of an $n$-nilpotent group and a rectangular band.
\end{thm}

\begin{proof}
($\rightarrow$) Let $\bS$ be an $n$-supernilpotent orthodox semigroup. Since an orthodox semigroup is regular, from Proposition \ref{ThmRegNilpSolv} it follows that $\bS$ is an $n$-supernilpotent completely simple semigroup. Since $\bS$ is orthodox and completely simple, Proposition \ref{HowieExer3.8} implies that it is isomorphic to a direct product $\bG\times(\mathbf{I}\times\mathbf{\Lambda})$ where $\mathbf{G}$ is a group and $(\mathbf{I}\times\mathbf{\Lambda})$ is a rectangular band. Now from Corollary \ref{SuperNilpOrthCSsg} it follows that $\bG$ is an $n$-nilpotent group. Hence, $\bS$ is isomorphic to a direct product of an $n$-nilpotent group and a rectangular band.

\noindent ($\leftarrow$) Let $\bS = \bG\times (\mathbf{I}\times\mathbf{\Lambda})$, where $\bG$ is an $n$-nilpotent
group and $\mathbf{I}\times\mathbf{\Lambda}$ is a rectangular band. From Proposition \ref{PropBandTimesGroup} it follows that $\bS$ is a completely simple orthodox semigroup. Since the group $\bG$ is $n$-nilpotent, from Corollary \ref{SuperNilpOrthCSsg} it follows that the semigroup $\bS$ is an $n$-supernilpotent semigroup. Hence, $\bS$ is an $n$-supernilpotent orthodox semigroup.
 \end{proof}

\begin{pro}\label{ThmOrthodoxNilpSolv}
Let $n\in\N$.
\begin{enumerate}
\item[(i)]  An orthodox semigroup is $n$-nilpotent if and only if it is isomorphic to a direct product of an $n$-nilpotent group and a rectangular band.
\item[(ii)] An orthodox semigroup is $n$-solvable if and only if it is isomorphic to a direct product of an $n$-solvable group and a rectangular band.
\end{enumerate}
\end{pro}

\begin{proof}
The proof is analogous to the proof of Theorem \ref{SupernilpOrthodox}.
\end{proof}

\begin{cor}
Let $n\in\N$. Then an orthodox semigroup is $n$-supernilpotent if and only if it is $n$-nilpotent.
\end{cor}

\begin{proof}
Let $n\in\N$. From Proposition \ref{ThmOrthodoxNilpSolv}(i) orthodox semigroup is $n$-nilpotent if and only if it is isomorphic to a direct product of $n$-nilpotent group and a rectangular band. Further by Proposition \ref{GroupNilpSupernilp} we know that it is equivalent to be isomorphic with $n$-supernilpotent group and a rectangular band and by Proposition \ref{GIL}(3) we know that it is equivalent for our orthodox semigroup to be $n$-supernilpotent.
\end{proof}

Inverse semigroups are regular semigroups such that each element has the unique inverse. They are obviously orthodox semigroups. The next statement has been proved independently by Kinyon and Stanovsk\'y in \cite{KD22}. 

\begin{pro}\label{inverse}
Let $n\in\N$. An inverse semigroup is $n$-supernilpotent if and only if it is isomorphic to a $n$-nilpotent group.
\end{pro}

\begin{proof}
Since inverse semigroup is regular then $n$-supernilpotent inverse semigroup is $n$-supernilpotent completely simple semigroup by Proposition \ref{ThmRegNilpSolv}(iii). Using that inverse semigroups are also orthodox we obtain that they are isomorphic to $\mathbf{G}\times(\mathbf{I}\times\mathbf{\Lambda})$ by Proposition \ref{HowieExer3.8}, where $\bG$ is a group, $\mathbf I$ is a left-zero semigroup and $\mathbf\Lambda$ is a right-zero semigroup. We know that all idempotents in inverse semigroups commute and idempotents in $\mathbf{G}\times(\mathbf{I}\times\mathbf{\Lambda})$
are of the form $(e,i,\lambda)$, where $e$ is the neutral element of the group $\bG$, $i\in I$ and $\lambda\in\Lambda$. Therefore, for all $j\in I$ and all $\kappa\in\Lambda$ should be $(e,i,\lambda)\cdot(e,j,\kappa)=(e,j,\kappa)\cdot(e,i,\lambda)$ and hence $(e,i,\kappa)=(e,j,\lambda)$ whence $i=j$ and
$\lambda=\kappa$. Therefore $\mathbf I$ and  $\mathbf\Lambda$ are trivial semigroups and
$\mathbf{G}\times(\mathbf{I}\times\mathbf{\Lambda})\cong\bG$. The opposite direction of the statement is obvious.
\end{proof}

\section{Acknowledgments}

We thank to P. Mayr and the referee for the useful discussion and
comments in preparation of this material.

\bibliographystyle{plain}

\end{document}